\documentclass{ifacconf}

\usepackage{graphicx}      
\usepackage{natbib}        
\usepackage{comment}
\usepackage{amsmath}
\usepackage{amssymb}
\usepackage{bm}
\usepackage{xcolor}
\usepackage{mathtools}
\usepackage{diffcoeff}[=v4]
\usepackage{empheq,etoolbox}

\usepackage{tikz}
\usetikzlibrary{shapes.geometric}

\newcommand{\bbR}{\mathbb{R}}
\renewcommand\d{\ensuremath{\mathrm{d}}}
\newcommand{\inpr}[3][]{\ensuremath{( #2, \, #3 )_{#1}}}

\DeclareMathOperator*{\grad}{grad}
\renewcommand{\div}{\operatorname{div}}
\DeclareMathOperator*{\curl}{curl}

\newtheorem{proposition}{Proposition}

\newtheorem{assumption}{Assumption}
\newtheorem{proof}{Proof}

\patchcmd{\subequations}
  {\theparentequation\alph{equation}}
  {\theparentequation.\alph{equation}}
  {}{}

\begin{document}
\begin{frontmatter}
\title{On the discrete equivalence of Lagrangian, Hamiltonian and mixed finite element formulations for linear wave phenomena} 

\author[ICA]{A. Brugnoli}
\author[TUB]{V. Mehrmann}

\address[ICA]{ICA, Université de Toulouse, ISAE-SUPAERO, MINES ALBI, UPS, INSA, CNRS, Toulouse, France}
\address[TUB]{Institut f\"ur Mathematik MA 4-5, Technische Universität Berlin, D-10623 Berlin, Germany}

\begin{abstract}
It is well known that the Lagrangian and Hamiltonian descriptions of field theories are equivalent at the discrete time level when variational integrators are used. Besides the symplectic Hamiltonian structure, many physical systems exhibit a Hamiltonian structure when written in mixed form. In this contribution, the discrete equivalence of Lagrangian, symplectic Hamiltonian and mixed formulations is investigated for linear wave propagation phenomena. Under compatibility conditions between the finite elements, the Lagrangian and mixed formulations are indeed equivalent. For the time discretization the leapfrog scheme and the implicit midpoint rule are considered. In mixed methods applied to wave problems the primal variable (e.g. the displacement in mechanics or the magnetic potential in electromagnetism) is not an unknown of the problem and is reconstructed a posteriori from its time derivative. When this reconstruction is performed via the trapezoidal rule, then these time-discretization methods lead to equivalent formulations. 
\end{abstract}

\begin{keyword}
Hamiltonian formulation, Lagrangian formulation, mixed finite elements. 
\end{keyword}

\end{frontmatter}

\section{Introduction}

Hamilton's principle of least action is the fundamental result behind classical field theories \cite{marsden2013introduction}. The equations may be obtained in Lagrangian or Hamiltonian form as the Euler-Lagrange equations are equivalent to the Hamiltonian equations. The latter equations are given by the canonical Poisson bracket but can be as well described by a Poisson bracket that defines a formally skew-adjoint differential operator. This is the viewpoint adopted in the seminal paper on port-Hamiltonian systems \cite{vanderSchaft2002} but also in mixed finite element formulation for dynamical systems (see for instance the seminal paper by \cite{geveci1988} on the velocity stress formulation of the wave equation). 

Recently in a series of papers, \cite{sanchez2017wave,sanchez2021elastodynamics,sanchez2022electro}, it was noticed that the continuous Galerkin and mixed finite element formulation preserve the Hamiltonian structure in general linear wave propagation phenomena. Discontinuous Galerkin method also preserve the Hamiltonain if the numerical fluxes are chosen in a suitable manner. Hybridizable discontinuous Galerkin methods reveal instead a dissipative Hamiltonian structure due to presence of stabilization terms. The Hamiltonian structure arises naturally from the Poisson brackets given by the weak formulation. This is perhaps not surprising as the  Euler-Lagrange equations lead to the strong and weak form of the equations of motion. Weak formulations based on finite elements naturally replicate this variational structure.

The equivalence of continuous Galerkin (i.e. Lagrangian) and mixed finite element discretization  has been thoroughly explored in \cite{joly2003}. Therein the authors mainly focused on the semi-discretization in space. By using variational time marching schemes the symplectic (or Poisson) structure can then be maintained at the fully discrete level. In this contribution we take the example of the Newmark time integration \cite{newmark1959method,kane2000newmark} applied to Lagrangian dynamics. In particular two instances of this integrator class will be considered: the leapfrog (or St\"ormer-Verlet) method and the implicit midpoint scheme. The first one is a partitioned Runge-Kutta method that is symplectic when applied to separable Hamiltonian systems. The second is the simplest Gauss Legendre collocation method that leads to exact energy conservation in the linear case \cite{kotyczka2019discrete,MehM19}. Mixed (or more generally port-Hamiltonian formulations) discard the primal variable of the problem (e.g. the displacement in elasticity or the magnetic potential in electromagnetism) when it does not contribute to the energy formulation. When this variable is reconstructed via the trapezoidal rule and the finite elements satisfy appropriate compatibility conditions, then the classical Hamiltonian formulation and the mixed formulation based on a Poisson bracket are equivalent.

The paper is organized as follows. In Sec. \ref{sec:cont_formulations} the wave and Maxwell equations are presented via the Hamilton principles. The associated Lagrangian, Hamiltonian and mixed formulations are detailed. The semi-discretization and the equivalence between continuous Galerkin and mixed formulation is discussed in Sec. \ref{sec:semidis_space}. In Sec. \ref{sec:time_disc} the fully discrete system is detailed for the wave equation only. The complete algebraic equivalence between Lagrangian, Hamiltonian and mixed scheme is investigated for the aforementioned time-integration schemes.

\paragraph*{Notation} We denote by $\Omega$ an open set of the space $\bbR^d$ of real $d$ dimensional vectors.  Given an Euclidean vector space $\mathbb{V}$, $||\cdot||$ denotes its Euclidean norm. Let $\mathbb T$ be a time interval. For a generic field $f : \Omega \times \mathbb T \rightarrow \mathbb{V}$, $||f||_{L^2}$ denotes its $L^2$ norm and  $\inpr[\Omega]{\cdot}{\cdot}$ the $L^2$ scalar product. The notation $\dot{q}$ is used both for the partial and ordinary time derivative of $q$. The notation $H(D) = \{u \in L^2 | Du \in L^2\}$ where $D\in \{\grad,\, \curl,\, \div \}$ denotes the standard Sobolev spaces. Whenever appropriate boundary conditions are considered we write $H_0(D)$. The notation $f^n:= f(n\Delta t)$ indicates evaluation of a function at a particular time  $n\Delta t$.

\section{Continuous formulations}\label{sec:cont_formulations}

In this section we detail the continuous formulation of the wave and Maxwell equations in order to show the affinities between the two. Indeed these models could be presented in greater generality using the language of exterior calculus and differential forms. We opted here for a  simpler representation based on  vector calculus formulations.

\subsection{Wave equation}
Consider $q$ to be the vertical displacement of an elastic membrane. The vibrations of the membrane are described by the linear wave equation, that corresponds to the \emph{Euler-Lagrange equations} associated with the 
\emph{Lagrangian}
\begin{equation*}
    L(q, \dot{q}) = \frac{1}{2} \int_\Omega \rho \dot{q}^2 -   k ||\nabla q||^2 \d\Omega,
\end{equation*}
where $\rho$ is the \emph{density} and $k$ is the \emph{stiffness}. Using \emph{Dirichlet boundary conditions}, this leads to the second order linear partial differential equation (PDE)
\begin{equation}\label{eq:wave_standard}
\begin{aligned}
    \rho \ddot{q} =  \nabla \cdot ( k \nabla q), \quad
    q|_{\partial\Omega} = 0.
\end{aligned}
\end{equation}
%

Denoting by $\delta_z$ the partial derivative with respect to the variable $z$, we consider the \emph{conjugate momentum} $p = \delta_{\dot{q}} L$ and the \emph{total energy or Hamiltonian}
\begin{equation}\label{Hamiltonian}
    H(q, p) =\int_\Omega p \dot{q} - L(q, p)\ \d\Omega.
\end{equation}
Note that the conjugated momentum and the Hamiltonian are formally defined using the \emph{Legendre transform} of the Lagrangian, see \cite{arnold2012geometrical}. 
The system can then be rewritten in \emph{symplectic Hamiltonian form} as
\begin{equation}\label{eq:wave_hamiltonian}
    \begin{bmatrix}
        \dot{p} \\
        \dot{q}
    \end{bmatrix} = \begin{bmatrix}
        0 & -I \\
        I & 0
    \end{bmatrix}
    \begin{bmatrix}
        \delta_{p} H \\
        \delta_{q} H
    \end{bmatrix},  \quad
    q|_{\partial\Omega} = 0,
\end{equation}
where $\delta_p H = \dot{q}$ and $\delta_{q} H = - \nabla \cdot (k\nabla q)$. This Hamiltonian system is also equivalent to a formulation using the velocity instead of the linear momentum
\begin{equation}\label{eq:wave_first}
    \begin{bmatrix}
        \rho & 0 \\
        0 & I
    \end{bmatrix}
    \begin{bmatrix}
        \dot{v} \\
        \dot{q}
    \end{bmatrix} = \begin{bmatrix}
        0 & -I \\
        I & 0
    \end{bmatrix}
    \begin{bmatrix}
        \dot{v} \\
         - \nabla \cdot (k\nabla q)
    \end{bmatrix},  \quad
    q|_{\partial\Omega} = 0.
\end{equation}
Another formulation, initially proposed in the finite element community \cite{geveci1988}, uses as variables the velocity $v = \dot{q}$ and the \emph{stress} $\bm{\sigma} = k \nabla q$, and leads to 
\begin{equation}\label{eq:wave_mixed}
    \begin{bmatrix}
    \rho & 0 \\
    0 & c \\
    \end{bmatrix}
    \begin{bmatrix}
        \dot{v}\\
        \dot{\bm\sigma}
    \end{bmatrix} = 
    \begin{bmatrix}
        0 & \nabla \cdot    \\
        \nabla & 0 \\
    \end{bmatrix}
    \begin{bmatrix}
        {v}\\
        {\bm\sigma}
    \end{bmatrix}, \qquad v|_{\partial\Omega}=0,
\end{equation}
where the \emph{compliance} $c:= k^{-1}$ has been introduced. In the finite element community this formulation is called a \emph{mixed formulation} and corresponds to a  Hamiltonian formulation, as the coefficient operator is skew-adjoint, see \cite{olver1993applications}. More generally this is also shown to be a \emph{port-Hamiltonian formulation}, see \cite{jacob2012linear}, as it has an underlying \emph{Dirac structure}. In this formulation, the Hamiltonian of the system is given by
\begin{equation}\label{Hamiltonianvsig}
    H({v}, \bm{\sigma}) = \frac{1}{2} \int_{\Omega} \rho v^2 + c||\bm{\sigma}||^2\ \d\Omega.
\end{equation}
From this velocity-stress formulation one can reduce the system to a second order formulation in the velocity variable only
\begin{equation}\label{eq:velocity_only}
    \rho \ddot{v} = \nabla \cdot (k\nabla v),
\end{equation}
or in the stress variable only
\begin{equation}\label{eq:stress_only}
    c\ddot{\bm\sigma} = \nabla ( \nu \nabla \cdot \bm\sigma),
\end{equation}
where $\nu$ is the \emph{specific volume} $\nu:=\rho^{-1}$. The velocity only formulation is just the time derivative of \eqref{eq:wave_standard} and it is therefore equivalent to it via integration in time. \\

A second mixed formulation is obtained by starting from the symplectic Hamiltonian formulation \eqref{eq:wave_first} and introducing the stress $\bm\sigma$ as an additional unknown via
\begin{equation}\label{eq:wave_mixed_ham}
\begin{aligned}
    \begin{bmatrix}
        \rho & 0 \\
        0 & I
    \end{bmatrix}
    \begin{bmatrix}
        \dot{v} \\
        \dot{q}
    \end{bmatrix} &= \begin{bmatrix}
        0 & -I \\
        I & 0
    \end{bmatrix}
    \begin{bmatrix}
        \dot{v} \\
         - \nabla \cdot \bm{\sigma}
    \end{bmatrix}, \\
    c \bm{\sigma} &= \nabla q.
\end{aligned}
\end{equation}
This is often the starting point for deriving \emph{mixed and hybridizable discontinuous finite element formulations} for the wave equation, see \cite{sanchez2017wave}.

It should be noted that all the 
discussed formulations have the same Hamiltonian, but are formulated in different variables. This may lead to advantages when considering limiting situations like letting the density go to zero of the stiffness go to $\infty$, see \cite{MehS23} for a  discussion in the finite dimensional case.

\subsection{Maxwell equations}
In the absence of electric charges, the Maxwell equations correspond to the Euler-Lagrange equations of the Lagrangian, see \cite{marsden2013introduction},
\begin{equation}
    L(\bm{A}, \dot{\bm{A}}) = \frac{1}{2} \int_\Omega \varepsilon ||\dot{\bm{A}}||^2 -   \mu^{-1} ||\nabla \times \bm{A}||^2\ \d\Omega,
\end{equation}
where $\varepsilon$ is the \emph{electric permittivity} and $\mu$ is the \emph{magnetic permeability}. This leads to the second order PDE
\begin{equation}\label{eq:max_standard}
\begin{aligned}
    \varepsilon \ddot{\bm{A}} =  -\nabla \times (\mu^{-1} \nabla \times \bm{A}), \qquad
    \bm{A} \times \bm{n}|_{\partial\Omega} = 0.
\end{aligned}
\end{equation}
Considering the \emph{conjugate momentum} $\bm{Y} = \delta_{\dot{\bm{A}}} L = \varepsilon\dot{\bm{A}} = -\bm{D}$ (where $\bm{D}$ is the \emph{electric flux density}), the total energy (Hamiltonian) has the form
\begin{equation}
    H(\bm{Y}, \bm{A}) =\int_\Omega \bm{Y} \dot{\bm{A}} - L(\bm{A}, \bm{Y})\ \d{\Omega}.
\end{equation}
The equation can then be rewritten in symplectic Hamiltonian form as
\begin{equation}\label{eq:max_hamiltonian}
    \begin{bmatrix}
        \dot{\bm{Y}} \\
        \dot{\bm{A}}
    \end{bmatrix} = \begin{bmatrix}
        0 & -I \\
        I & 0
    \end{bmatrix}
    \begin{bmatrix}
        \delta_{\bm{Y}} H \\
        \delta_{\bm{A}} H
    \end{bmatrix},
\end{equation}
where $\delta_{\bm{Y}} H = \dot{\bm{A}}$ and $\delta_{\bm{A}} H = \nabla \times (\mu^{-1}\nabla \times \bm{A})$. This Hamiltonian system is also equivalent to a formulation using the \emph{electric field} $\bm{E} = - \varepsilon^{-1} \bm{Y}$,
\begin{equation}\label{eq:max_first}
    \begin{bmatrix}
        \varepsilon & 0 \\
        0 & I
    \end{bmatrix}
    \begin{bmatrix}
        \dot{\bm{E}} \\
        \dot{\bm{A}}
    \end{bmatrix} = \begin{bmatrix}
        0 & I \\
        -I & 0
    \end{bmatrix}
    \begin{bmatrix}
        \dot{\bm{E}} \\
         \nabla \times (\mu^{-1}\nabla \times \bm{A})
    \end{bmatrix}.
\end{equation}
The mixed formulation uses as variables the electric and \emph{magnetic field} $\bm{H} = \mu^{-1} \nabla \times \bm{A}$,
\begin{equation}\label{eq:max_mixed}
    \begin{bmatrix}
        \varepsilon & 0 \\
        0 & \mu
    \end{bmatrix}
    \begin{bmatrix}
        \dot{\bm{E}}\\
        \dot{\bm{H}}
    \end{bmatrix} = 
    \begin{bmatrix}
        0 & \nabla \times    \\
        -\nabla \times & 0 \\
    \end{bmatrix}
    \begin{bmatrix}
        {\bm{E}}\\
        {\bm{H}}
    \end{bmatrix}, \qquad \bm{E} \times \bm{n}|_{\partial\Omega}=0.
\end{equation}
The Hamiltonian of this system is given by
\begin{equation*}
    H = \frac{1}{2} \int_{\Omega} \varepsilon ||\bm{E}||^2 + \mu||\bm{H}||^2\ \d\Omega.
\end{equation*}
From the mixed formulation formulation one can  reduce the system to one for the electric field only
\begin{equation}\label{eq:electric_only}
    \varepsilon \ddot{\bm{E}} = - \nabla \times (\mu^{-1}\nabla \times \bm{E}),
\end{equation}
or for the magnetic field only 
\begin{equation}\label{eq:magnetic_only}
    \mu \ddot{\bm{H}} = - \nabla \times ( \varepsilon^{-1} \nabla \times \bm{H}).
\end{equation}
An alternative formulation is obtained using the variables $\bm{E}, \bm{A}$  but introducing the magnetic field definition
\begin{equation}\label{eq:max_mixed_ham}
\begin{aligned}
    \begin{bmatrix}
        \varepsilon & 0 \\
        0 & I
    \end{bmatrix}
    \begin{bmatrix}
        \dot{\bm{E}} \\
        \dot{\bm{A}}
    \end{bmatrix} &= \begin{bmatrix}
        0 & I \\
        -I & 0
    \end{bmatrix}
    \begin{bmatrix}
        \dot{\bm{E}} \\
         \nabla \times \bm{H}
    \end{bmatrix}, \\
    \mu \bm{H} &= \nabla \times \bm{A}.
\end{aligned}
\end{equation}
This is the formulation used to devise the \emph{hybridizable discontinuous Galerkin method} in \cite{sanchez2022electro}.

\section{Semi-discretization in space}\label{sec:semidis_space}

In this section we consider the semi-discretization of the wave and Maxwell equations. Once again the discussion in more or less analogous in the two cases and may be unified via the formalism of Finite Elements Exterior calculus \cite{arnold2006acta}. 

\subsection{The semi-discrete wave equation}
The classical Lagrangian formulation of the wave equation can be discretized using conforming finite element spaces $V_{h, 0}(\grad) \subset H_0(\grad)$, see \cite{arnold2006acta}, leading to the following weak formulation: 

Find $q_h \in V_{h, 0}(\grad)$ such that 
\begin{equation}\label{eq:wave_standard_weak}
    \inpr[\Omega]{\psi_h}{\rho\ddot{q}_h} = -\inpr[\Omega]{\nabla \psi_h}{k \nabla q_h}\ \mbox{\rm for all}\ \psi_h \in V_{h, 0}(\grad).
\end{equation}

To establish the various equivalences between the different formulations we make the following assumption.
\begin{assumption}\label{ass1}
    All the physical coefficients are assumed to be constant.
\end{assumption}
Given Asumption~\ref{ass1}, the weak form of the Hamiltonian formulation \eqref{eq:wave_first} is clearly equivalent to \eqref{eq:wave_standard_weak} when one picks $q_h, v_h \in V_{h, 0}(\grad)$. For  the mixed formulations, two possibilities via integration by parts can be pursued.
Integration by parts on the first line leads to the formulation: 

Find $(v_h, \bm{\sigma}_h) \in V_{h, 0}(\grad) \times \bm{W}_h \subset H_0(\grad) \times L^2$ such that
\begin{equation}\label{eq:wave_mixed_grad}
\begin{aligned}
    \inpr[\Omega]{\psi_h}{\rho\dot{v}_h} &= -\inpr[\Omega]{\nabla \psi_h}{\bm{\sigma}_h}, \\
    \inpr[\Omega]{\bm{\xi}_h}{c\dot{\bm{\sigma}}} &=  \inpr[\Omega]{\bm{\xi}_h}{\nabla v_h},
\end{aligned} \
\begin{aligned}
    &\ \mbox{\rm for all}\ \psi_h \in V_{h, 0}(\grad), \\
    &\ \mbox{\rm for all}\ \bm{\xi}_h \in \bm{W}_h.
\end{aligned}
\end{equation}
If the space $\bm{W}_h$ is such that $\bm{W}_h \subset \nabla V_{h, 0}(\grad)$, and since the coefficients are assumed to be constants, the second equation holds pointwise, i.e., $c\dot{\bm{\sigma}} = \nabla v_h$. One can then take the time derivative of the first equation and use the second equation to obtain the weak formulation: 

Find $v_h \in V_{h, 0}(\grad)$
\begin{equation}\label{eq:velocity_only_weak}
    \inpr[\Omega]{\psi_h}{\rho \ddot{v}} = -\inpr[\Omega]{\nabla \psi_h}{\nabla v_h}, \ \mbox{\rm for all}\ \psi_h \in V_{h, 0}(\grad).
\end{equation}
This is exactly the continuous Galerkin weak form of the second order formulation for the velocity and it is therefore equivalent to the classical Lagrangian weak formulation \eqref{eq:wave_standard_weak} when its time integration is performed.

The second possibility consists in integrating by parts the second line, leading to the weak formulation: 

Find $(v_h, \bm{\sigma}_h) \in W_h \times \bm{V}_h(\div) \subset L^2 \times H(\div)$ such that  
\begin{equation}\label{eq:wave_mixed_div}
\begin{aligned}
    \inpr[\Omega]{\psi_h}{\rho\dot{v}_h} &= \inpr[\Omega]{\psi_h}{\nabla \cdot \bm{\sigma}_h}, \\
    \inpr[\Omega]{\bm{\xi}_h}{c\dot{\bm{\sigma}}} &=  -\inpr[\Omega]{\nabla \cdot \bm{\xi}_h}{v_h},
\end{aligned} \
\begin{aligned}
    &\ \mbox{\rm for all}\ \psi_h \in W_h, \\
    &\ \mbox{\rm for all}\ \bm{\xi}_h \in \bm{V}_h(\div).
\end{aligned}
\end{equation}
Analogous to the previous case, if the  space $W_h$ satisfies ${W}_h \subset \nabla \cdot \bm{V}_h(\div)$, then the first equation holds pointwise, i.e., $\rho \dot{v}_h = \nabla \cdot \bm{\sigma}_h$, and then taking the time derivative of the second equation one obtains
 \begin{equation}\label{eq:stress_only_weak}
     \inpr[\Omega]{\bm{\xi}_h}{c \ddot{\bm{\sigma}}_h} = - \inpr[\Omega]{\nabla \cdot \bm{\xi}_h}{\nu \nabla \cdot \bm{\sigma}_h}.
 \end{equation}
An alternative mixed formulation is obtained using \eqref{eq:wave_mixed_ham},
\begin{equation}\label{eq:wave_first_mixed_div}
    \begin{aligned}
        \begin{aligned}
    \inpr[\Omega]{\psi_h}{\rho\dot{v}_h} &= \inpr[\Omega]{\psi_h}{\nabla \cdot \bm{\sigma}_h}, \\
    \inpr[\Omega]{\psi_h}{\dot{q}_h} &=  \inpr[\Omega]{\psi_h}{v_h}, \\
    \inpr[\Omega]{\bm{\xi}_h}{c{\bm{\sigma}_h}} &=  -\inpr[\Omega]{\nabla \cdot \bm{\xi}_h}{q_h},
\end{aligned} \
\begin{aligned}
    &\ \mbox{\rm for all}\ \psi_h \in W_h, \\
    &\ \mbox{\rm for all}\ \psi_h \in W_h, \\
    &\ \mbox{\rm for all}\ \bm{\xi}_h \in \bm{V}_h(\div).
\end{aligned}
    \end{aligned}
\end{equation}
This formulation is equivalent to \eqref{eq:wave_mixed_div} by taking the time derivative of the last equation.
 
The equivalence between different formulations is discussed in  detail in \cite{joly2003}. Therein, however, the geometric interpretation of the different formulations is not discussed. 
Geometry plays, however, a fundamental role as the connection between  \eqref{eq:wave_mixed_grad} and \eqref{eq:wave_mixed_div} is given by the \emph{Hodge operator}, that maps differential forms to their dual space isomorphically. In a classical mixed finite element formulation, the Hodge operator is expressed by a projection between dual spaces of finite elements.  This projection entails a loss of information. An isomorphic Hodge star requires dual meshes and this is the approach in the Discrete Exterior Calculus, \cite{hirani2003discrete}.

 In the following  we make a second assumption.

 \begin{assumption}\label{ass2}
The finite element spaces satisfy the compatibility conditions
 \begin{equation*}
     \bm{W}_h \subset \nabla V_{h, 0}(\grad), \qquad W_h \subset \nabla \cdot \bm{V}_h(\div).
 \end{equation*}
 \end{assumption}

 \subsection{The semi-discrete Maxwell equations}
The Lagrangian formulation of the Mawxell equations can be discretized using conforming finite elements $\bm{V}_{h, 0}(\curl) \subset H_0(\curl)$, leading to the following weak formulation: 

Find $\bm{A}_h \in V_{h, 0}(\curl)$ such that for all $ \bm{\psi}_h \in \bm{V}_{h, 0}(\curl)$,
\begin{equation}\label{eq:max_standard_weak}
    \inpr[\Omega]{\bm{\psi}_h}{\varepsilon\ddot{\bm{A}}_h} = -\inpr[\Omega]{\nabla \times \bm{\psi}_h}{\mu^{-1}\nabla \times \bm{A}}.
\end{equation}
The weak form of the Hamiltonian formulation \eqref{eq:max_hamiltonian} and \eqref{eq:max_first} is clearly equivalent to \eqref{eq:max_standard_weak} when one picks $\bm{E}_h^0, \bm{A}_h^0 \in \bm{V}_{h, 0}(\curl)$. The integration by parts on the first line leads to the formulation: 

Find $(\bm{E}_h, \bm{A}_h) \in \bm{V}_{h, 0}(\curl) \times \bm{W}_h \subset H_0(\curl) \times L^2(\bbR^d)$ such that
\begin{equation}\label{eq:max_mixed_first}
\begin{aligned}
    \inpr[\Omega]{\bm{\psi}_h}{\varepsilon\dot{\bm{E}}_h} &= \inpr[\Omega]{\nabla \times \bm{\psi}_h}{\bm{H}_h}, \\
    \inpr[\Omega]{\bm{\xi}_h}{\mu\dot{\bm{H}}} &=  -\inpr[\Omega]{\bm{\xi}_h}{\nabla \times \bm{E}_h},
\end{aligned} \
\begin{aligned}
    &\mbox{\rm for all}\ \bm{\psi}_h \in \bm{V}_{h, 0}(\curl), \\
    &\mbox{\rm for all}\ \bm{\xi}_h \in \bm{W}_h.
\end{aligned}
\end{equation}
If the  space $\bm{W}_h$ is such that $\bm{W}_h \subset \nabla \times V_{h, 0}(\curl)$, then the second equation holds pointwise, i.e., $\mu \dot{\bm{H}} = -\nabla \times \bm{E}_h$. The weak formulation when using only the electric field takes the form: 

Find $\bm{E}_h \in V_{h, 0}(\curl)$ such that for all 
$ \bm{\psi}_h \in \bm{V}_{h, 0}(\curl)$
\begin{equation}\label{eq:electric_only_weak}
    \inpr[\Omega]{\bm{\psi}_h}{\varepsilon\ddot{\bm{E}}_h} = -\inpr[\Omega]{\nabla \times \bm{\psi}_h}{\mu^{-1}\nabla \times \bm{E}}.
\end{equation}
This is 
equivalent to the classical Lagrangian weak formulation \eqref{eq:wave_standard_weak}.

The second possibility consists in integrating by parts the second line, leading to the weak formulation: 

Find $(\bm{E}_h, \bm{H}_h) \in \bm{W}_h \times \bm{V}_h(\curl) \subset L^2 \times H(\curl)$ such that
\begin{equation}\label{eq:max_mixed_second}
\begin{aligned}
    \inpr[\Omega]{\bm{\psi}_h}{\varepsilon\dot{\bm{E}}_h} &= \inpr[\Omega]{\bm{\psi}_h}{\nabla \times  \bm{H}_h}, \\
    \inpr[\Omega]{\bm{\xi}_h}{\mu\dot{\bm{H}}} &=  -\inpr[\Omega]{\nabla \times \bm{\xi}_h}{ \bm{E}_h},
\end{aligned} \
\begin{aligned}
    &\mbox{\rm for all}\ \bm{\xi}_h \in \bm{W}_h, \\
    &\mbox{\rm for all}\ \bm{\psi}_h \in \bm{V}_{h, 0}(\curl). 
\end{aligned}
\end{equation}
If $\bm{W}_h \subset \nabla \times \bm{V}_h(\div)$, then the first equation holds pointwise, i.e., $\varepsilon \dot{\bm{E}}_h = \nabla \times \bm{H}_h$. Then taking the time derivative of the second equation one obtains
 \begin{equation}\label{eq:magnetic_only_weak}
     \inpr[\Omega]{\bm{\xi}_h}{\mu \ddot{\bm{H}}_h} = - \inpr[\Omega]{\nabla \times \bm{\xi}_h}{\varepsilon^{-1} \nabla \cdot \bm{\sigma}_h}.
 \end{equation}
 An alternative mixed formulation is obtained using formulation,
\begin{equation}\label{eq:max_mixed_ham_second}
\begin{aligned}
    \inpr[\Omega]{\bm{\psi}_h}{\dot{\bm{E}}_h} &= \inpr[\Omega]{\bm{\psi}_h}{\nabla \times{\bm{H}}_h}, \\
    \inpr[\Omega]{\bm{\psi}_h}{\dot{\bm{A}}_h} &= \inpr[\Omega]{\bm{\psi}_h}{{\bm{E}_h}}, \\
    \inpr[\Omega]{\bm{\xi}_h}{\mu \bm{H}_h} &= \inpr[\Omega]{\bm{\xi}_h}{\nabla \times \bm{A}_h}.
\end{aligned}
\end{equation}
 The same considerations concerning the time discretization carry over to the Maxwell equations, leading to a complete equivalence of discrete-time formulations of Lagrangian, Hamiltonian and mixed finite element descriptions for electromagnetic phenomena.

\section{Time discretization}\label{sec:time_disc}

For the time discretization of the different second order formulations, two different schemes will be considered: the \emph{symplectic leapfrog method} and the \emph{implicit midpoint rule}. These are particular instances of the \emph{Newmark method}, originally developed for structural dynamics in \cite{newmark1959method}. A general Newmark scheme applied to a second order system has the form
\begin{subequations}\label{eq:newmark}
\begin{empheq}{align}
    \inpr[\Omega]{\psi_h}{\rho a_h^{n+1}} &= -\inpr[\Omega]{\nabla \psi_h}{\nabla q_h^{n+1}}, \label{eq:newmark_1} \\
        \frac{v_h^{n+1} - v_h^n}{\Delta t} &=  \gamma a_h^{n+1} + (1 - \gamma) a_h^{n}, \label{eq:newmark_2} \\
        \frac{q_h^{n+1} - q_h^n}{\Delta t}&= v_h^n + \Delta t(\beta a_h^{n+1} + (\frac{1}{2}-\beta) a_h^n), \label{eq:newmark_3}
\end{empheq}
\end{subequations}
where $a_h^n$ denotes the acceleration at time $t^n$. For ease of presentation the results will be presented for the wave equation only as everything carries over to the Maxwell equations case (and \textit{mutatis mutandis} to the linear elastodynamics problem and derived models, like beams and plates structural models).

\subsection{The leapfrog scheme}
The Newmark scheme is equivalent to the leapfrog scheme when $\gamma=\frac{1}{2}$ and $\beta=0$, leading to the following system
\begin{equation}\label{eq:wave_standard_sv}
    \inpr[\Omega]{\psi_h}{\rho(q^{n+1}_h - 2q^n_h + q^{n-1}_h)} = -\Delta t^2 \inpr[\Omega]{\nabla \psi_h}{k \nabla q_h^n}.
\end{equation}
This is shown by considering two consecutive updates for the displacement \eqref{eq:newmark_3} and using \eqref{eq:newmark_2}. This scheme is also equivalent to the \emph{St\"ormer-Verlet method} applied to the weak formulation of \eqref{eq:wave_first} (see e.g. \cite{hairer2003})
\begin{subequations}\label{eq:wave_first_sv}
\begin{empheq}{align}
    \inpr[\Omega]{\psi_h}{\rho(v_{h, L}^{n+\frac{1}{2}} - v_{h, L}^{n-\frac{1}{2}})} &= -\Delta t \inpr[\Omega]{\nabla \psi_h}{k \nabla q_{h, L}^n}, \label{eq:wave_first_sv_1}\\
    q_{h, L}^{n+1} - q_{h, L}^{n} &= \Delta t v_{h, L}^{n+\frac{1}{2}},
\end{empheq}
\end{subequations}
where the subscript $L$ stands for Lagrangian description. For the mixed formulation \eqref{eq:wave_mixed_grad} one obtains
\begin{subequations}\label{eq:wave_mixed_grad_sv}
\begin{empheq}{align}
\inpr[\Omega]{\psi_h}{\rho(v_{h, M}^{n+\frac{1}{2}} - v_{h, M}^{n-\frac{1}{2}})} &= -\Delta t \inpr[\Omega]{\nabla \psi_h}{\bm{\sigma}_{h, M}^{n}}, \label{eq:wave_mixed_grad_sv_1}\\
\inpr[\Omega]{\bm{\xi}_h}{c(\bm{\sigma}_{h, M}^{n+1} - \bm{\sigma}_{h, M}^n)} &=  \Delta t \inpr[\Omega]{\bm{\xi}_h}{\nabla v_{h, M}^{n+\frac{1}{2}}}, \label{eq:wave_mixed_grad_sv_2}
\end{empheq}
\end{subequations}
where the subscript $M$ stands for mixed formulation.
\begin{proposition}\label{pr:equivalence_sv}
 Suppose that $\sigma_{h, M}^0 = k \nabla q_{h, L}^0$ and that the field $q_{h, M}^{n+1}$ with $q_{h, M}^{0}= q_{h, L}^{0}$ is reconstructed via the \emph{trapezoidal rule}
 \[
 q_{h, M}^{n+1} = q_{h, M}^{n} + \frac{\Delta t}{2} (v_{h, M}^n + v_{h, M}^{n+1}),
 \]
 then the formulations \eqref{eq:wave_first_sv} and \eqref{eq:wave_mixed_grad_sv} are equivalent.
\end{proposition}
\begin{proof}
It is sufficient to show that $\bm{\sigma}_{h, M}^n=\nabla q_{h, L}^n, \; \mbox{\rm for all}\ n$ as this implies that the dynamic equations \eqref{eq:wave_first_sv_1} and \eqref{eq:wave_mixed_grad_sv_1} are the same. The reconstruction of $q_h$ is the same by assumption.

 Using the trapezoidal rule for $q$ and \eqref{eq:wave_mixed_grad_sv_1}  one has 
\[
\bm{\sigma}_{h, M}^{n+1} - \bm{\sigma}_{h, M}^n = k \nabla (q_{h, M}^{n+1} - q_{h, M}^{n}).
\]
Since $\sigma_{h, M}^0 = k \nabla q_{h, L}^0$, the result is obtained by recursion.
\end{proof}

Formulation \eqref{eq:wave_mixed_grad_sv} is also equivalent to a staggered leapfrog discretization of \eqref{eq:velocity_only_weak} given by
\begin{equation}
    \inpr[\Omega]{\psi_h}{\rho(v_h^{n+\frac{3}{2}} - 2 v_h^{n+\frac{1}{2}} + v_h^{n-\frac{1}{2}})} = -\Delta t^2 \inpr[\Omega]{\nabla \psi_h}{k \nabla v_h^{n+\frac{1}{2}}}.
\end{equation}
To see this, it is sufficient to take the difference between two consecutive step for the velocity update \eqref{eq:wave_mixed_grad_sv_1} and use \eqref{eq:wave_mixed_grad_sv_2}. In a dual manner, the application of the St\"ormer-Verlet scheme to the mixed formulation \eqref{eq:wave_mixed_div}
\begin{equation}\label{eq:wave_mixed_div_sv}
\begin{aligned}
    \inpr[\Omega]{\psi_h}{\rho(v_{h, M}^{n+\frac{1}{2}} - v_{h, M}^{n-\frac{1}{2}})} &= \Delta t\inpr[\Omega]{\psi_h}{\nabla \cdot \bm{\sigma}_h^n}, \\
    \inpr[\Omega]{\bm{\xi}_h}{c(\bm{\sigma}_{h, M}^{n+1} - \bm{\sigma}_{h, M}^n)} &=  - \Delta t\inpr[\Omega]{\nabla \cdot \bm{\xi}_h}{v_h^{n+\frac{1}{2}}},
\end{aligned}
\end{equation}
is equivalent to the leapfrog scheme applied to \eqref{eq:stress_only_weak}
 \begin{equation}\label{eq:stress_only_weak_sv}
     \inpr[\Omega]{\bm{\xi}_h}{c(\bm{\sigma}_h^{n+1} - 2 \bm{\sigma}_h^n + \bm{\sigma}_h^{n-1})} = - \inpr[\Omega]{\nabla \cdot \bm{\xi}_h}{\nu \nabla \cdot \bm{\sigma}_h^n}.
 \end{equation}
The application of the St\"ormer-Verlet to the alternative mixed formulation \eqref{eq:wave_first_mixed_div} leads to 
\begin{equation}\label{eq:wave_first_mixed_div_sv}
        \begin{aligned}
    \inpr[\Omega]{\psi_h}{\rho({v}_h^{n+\frac{1}{2}} - {v}_h^{n-\frac{1}{2}})} &= \Delta t\inpr[\Omega]{\psi_h}{\nabla \cdot \bm{\sigma}_h^n}, \\
    \inpr[\Omega]{\psi_h}{q_h^{n+1} - q_h^{n}} &=  \Delta t\inpr[\Omega]{\psi_h}{v_h^{n+\frac{1}{2}}}, \\
    \inpr[\Omega]{\bm{\xi}_h}{c{\bm{\sigma}_h^{n+1}}} &=  -\inpr[\Omega]{\nabla \cdot \bm{\xi}_h}{q_h^{n+1}}.
\end{aligned} 
\end{equation}
Again this scheme is equivalent to \eqref{eq:wave_mixed_div_sv}, when $q_h$ is obtained from $v_h$ using the trapezoidal rule.

\subsection{The implicit midpoint rule}
For $\gamma=\frac{1}{2}$ and $\beta=\frac{1}{4}$ the Newmark scheme  leads to the implicit midpoint rule 
\begin{equation}
    \frac{v_{h, L}^{n+1} - v_{h, L}^n}{\Delta t} = a_h^{n+\frac{1}{2}}, \qquad \frac{q_{h, L}^{n+1} - q_{h, L}^n}{\Delta t}  =  v_h^{n + \frac{1}{2}},
\end{equation}
 where the notation $f_h^{n+\frac{1}{2}} := \frac{f_h^{n+1} + f_h^n}{2}$ has been used. Then system \eqref{eq:newmark} is rewritten as 
\begin{subequations}\label{eq:wave_first_im}
\begin{empheq}{align}
    \inpr[\Omega]{\psi_h}{\rho (v_{h, L}^{n+1} - v_{h, L}^n)} &= -\Delta t \inpr[\Omega]{\nabla \psi_h}{\nabla q_{h, L}^{n+\frac{1}{2}}}, \label{eq:wave_first_im_1}\\
    q_{h, L}^{n+1} - q_{h, L}^n &=   \Delta t v_h^{n + \frac{1}{2}}. \label{eq:wave_first_im_2}
\end{empheq}
\end{subequations}
The midpoint rule applied to the mixed discretization \eqref{eq:wave_mixed_grad} leads to 
\begin{subequations}\label{eq:wave_mixed_grad_im}
\begin{empheq}{align}
    \inpr[\Omega]{\psi_h}{\rho (v_{h, M}^{n+1} - v_{h, M}^n)} &= -\Delta t \inpr[\Omega]{\nabla \psi_h}{\bm{\sigma}_{h, M}^{n+\frac{1}{2}}}, \label{eq:wave_mixed_grad_im_1} \\
    \inpr[\Omega]{\bm{\xi}_h}{c(\bm{\sigma}_{h, M}^{n+1} - \bm{\sigma}_{h, M}^n)} &=  \Delta t \inpr[\Omega]{\bm{\xi}_h}{\nabla v_{h, M}^{n+\frac{1}{2}}}.\label{eq:wave_mixed_grad_im_2}
\end{empheq}
\end{subequations}
\begin{proposition}\label{pr:equivalence_im}
Under the assumptions of Proposition \ref{pr:equivalence_sv}  the formulations \eqref{eq:wave_first_im} and \eqref{eq:wave_mixed_grad_im} are equivalent.
\end{proposition}
\begin{proof}
The proof is analogous to that of Proposition \ref{pr:equivalence_sv}.
\end{proof}
The implicit midpoint rule applied to \eqref{eq:wave_mixed_grad_im} is equivalent to the following iteration for \eqref{eq:velocity_only_weak}
\begin{equation}
    \inpr[\Omega]{\psi_h}{\rho(v_h^{n+1} - 2 v_h^{n} + v_h^{n-1})} = -\Delta t^2 \inpr[\Omega]{\nabla \psi_h}{k \nabla \widehat{v}_h^{n}},
\end{equation}
where 
$$\widehat{v}_h^{n}:= \frac{1}{4}(v_h^{n+1} + 2 v_h^{n} + v_h^{n-1}).$$
In an analogous fashion as for the leapfrog scheme, this is shown by taking two consecutive time steps for the velocity update \eqref{eq:wave_first_im_1} and using \eqref{eq:wave_first_im_2}. For the dual formulation, the implicit midpoint rule
\begin{equation}\label{eq:wave_mixed_div_im}
\begin{aligned}
    \inpr[\Omega]{\psi_h}{\rho(v_{h, M}^{n+1} - v_{h, M}^{n})} &= \Delta t\inpr[\Omega]{\psi_h}{\nabla \cdot \bm{\sigma}_h^{n+\frac{1}{2}}}, \\
    \inpr[\Omega]{\bm{\xi}_h}{c(\bm{\sigma}_{h, M}^{n+1} - \bm{\sigma}_{h, M}^n)} &=  - \Delta t\inpr[\Omega]{\nabla \cdot \bm{\xi}_h}{v_h^{n+\frac{1}{2}}},
\end{aligned}
\end{equation}
leads to the following update when applied to \eqref{eq:stress_only_weak}
 \begin{equation}
     \inpr[\Omega]{\bm{\xi}_h}{c(\bm{\sigma}_h^{n+1} - 2 \bm{\sigma}_h^n + \bm{\sigma}_h^{n-1})} = - \inpr[\Omega]{\nabla \cdot \bm{\xi}_h}{\nu \nabla \cdot \widehat{\bm{\sigma}}_h^n},
 \end{equation}
with 
\[ \widehat{\bm{\sigma}}_h^n:= \frac{1}{4}(\bm{\sigma}_h^{n+1} + 2 \bm{\sigma}_h^n + \bm{\sigma}_h^{n-1}).\]

The application of the implicit midpoint rule to system \eqref{eq:wave_first_mixed_div} leads to 
\begin{equation}\label{eq:wave_first_mixed_div_im}
\begin{aligned}
    \inpr[\Omega]{\psi_h}{\rho({v}_h^{n+1} - {v}_h^{n})} &= \Delta t\inpr[\Omega]{\psi_h}{\nabla \cdot \bm{\sigma}_h^{n+\frac{1}{2}}}, \\
    \inpr[\Omega]{\psi_h}{q_h^{n+1} - q_h^{n}} &=  \Delta t\inpr[\Omega]{\psi_h}{v_h^{n+\frac{1}{2}}}, \\
    \inpr[\Omega]{\bm{\xi}_h}{c{\bm{\sigma}_h^{n+\frac{1}{2}}}} &=  -\inpr[\Omega]{\nabla \cdot \bm{\xi}_h}{q_h^{n+\frac{1}{2}}},
\end{aligned} 
\end{equation}
which is equivalent to \eqref{eq:wave_mixed_div_im} when $q_h$ is reconstructed via the trapezoidal rule.

 \section{Conclusion}
 Several different formulations of standard and mixed finite element discretizations as well as appropriate time discretization methods for linear wave phenomena have been compared. It is shown that with appropriate choices of space discretization methods, Hamiltonian, Lagrangian and mixed formulations lead to equivalent formulations  and also appropriate time discretization schemes lead to equivalent schemes. 

A natural question that arises is whether the discrete equivalence carries over to the numerical linear algebra level. Indeed if one uses the Newmark integrator on a second order system, a positive definite system has to be solved. On the other hand the application of the implicit midpoint to a mixed formulation leads to positive mass matrix perturbed by a small skew-symmetric matrix. Is seems then possible to implement equivalent iterative solvers for these two different problems. Indeed it was shown in \cite{guducu2022} that some iterative schemes for the Hamiltonian formulation indeed lead to similar convergence rate as the conjugate gradient (applied to symmetric positive definite problems).

\begin{ack}
The first author would like to thank Enrico Zampa from Trento University for insightful discussions on the topic.
\end{ack}

\bibliography{biblio}

\end{document}